\documentclass[a4paper,11pt]{article}
\usepackage{amsmath,amssymb,amsfonts} 
\usepackage[all]{xy}
\DeclareMathOperator{\rank}{rank}
\renewcommand{\AA}{\mathbb{A}}
\newcommand{\CC}{\mathbb{C}}

\newcommand{\PP}{\mathbb{P}}
\newcommand{\QQ}{\mathbb{Q}}
\newcommand{\ZZ}{\mathbb{Z}}

\newcommand{\Ce}{\mathcal{C}}
\newcommand{\Oh}{\mathcal{O}}

\newcommand{\al}{\alpha}
\newcommand{\be}{\beta}
\newcommand{\gam}{\gamma}

\newcommand{\fie}{\varphi}
\newcommand{\Fie}{\Phi}

\newcommand{\la}{\lambda}
\newcommand{\git}{/\!\!/}

\DeclareMathOperator{\Proj}{Proj}
\DeclareMathOperator{\Pic}{Pic}

\DeclareMathOperator{\Tors}{Tors}

\DeclareMathOperator{\Image}{Im}

\newtheorem{thm}{Theorem}[section]
\newtheorem{mainthm}[thm]{Main Theorem}
\newtheorem{lemma}[thm]{Lemma}
\newtheorem{prop}[thm]{Proposition}
\newtheorem{cor}[thm]{Corollary}

\newtheorem{rmk}[thm]{Remark}
\newenvironment{pf}{\paragraph{Proof}}{\par\medskip}
\newenvironment{pfcor}{\paragraph{Proof of corollary}}{\par\medskip}

\newenvironment{pfmainthm}{\paragraph{Proof of main theorem}}{\par\medskip}

\title{Extending hyperelliptic K3 surfaces,\\ and Godeaux surfaces with torsion $\ZZ/2$}
\author{Stephen Coughlan}
\date{}
\begin{document}
\maketitle
\begin{abstract}
We study the extension of a hyperelliptic K3 surface to a Fano $6$-fold. This determines a family of surfaces of general type with $p_g=1$, $K^2=2$ and hyperelliptic canonical curve, where each surface is a weighted complete intersection inside a Fano $6$-fold. Finally, we use these hyperelliptic surfaces to determine an $8$-parameter family of Godeaux surfaces with torsion $\ZZ/2$.
\end{abstract}

\section{Introduction}
We begin by studying extensions of certain hyperelliptic K3 surfaces. These surfaces are the hyperelliptic degeneration of the symmetric determinantal quartic surfaces studied in \cite{part1}. Let $D$ be a hyperelliptic curve of genus $3$ and $A$ be an ineffective theta characteristic on $D$. We study extensions of the graded ring \[R(D,A)=\bigoplus_{n\ge0}H^0(D,\Oh_D(nA))\] with $\Proj R(D,A)$ defining $D\subset\PP(2^3,3^4,4)$. In the first instance, there is an extension of $D$ to a hyperelliptic K3 surface $T\subset\PP(2^4,3^4,4)$ with $10\times\frac12(1,1)$ points and containing $D$ as a weighted hyperplane section of degree $2$. Sections \ref{sec!hellcurve} and \ref{sec!K3surface} of this paper treat graded rings over hyperelliptic curves and K3 surfaces respectively, by working relative to the hyperelliptic double covering.

Now, the K3 surface $T$ is the elephant hyperplane section of a Fano $3$-fold $W\subset\PP(1,2^4,3^4,4)$ with $10\times\frac12$ points. We think of this $3$-fold as an extension of $T$, and make further extensions up to a Fano $6$-fold $W^6$, with each successive $W_i$ containing $T$ as an appropriate number of hyperplane sections. This leads to the tower
\[D\subset T\subset W^3\subset W^4\subset W^5\subset W^6\subset\PP(1^4,2^4,3^4,4).\]
This is the hyperelliptic degeneration of the symmetric determinantal quartic extensions constructed in \cite{part1}. Thus we can hope for a one-to-one correspondence of moduli between the surface $T$ and the $6$-fold $W^6$. Indeed, we have
\begin{mainthm}\label{thm!hell} For each quasismooth hyperelliptic K3 surface $T\subset\PP(2^4,3^4,4)$ with $10\times\frac12$ points, there is a unique extension to a quasismooth Fano $6$-fold $W\subset\PP(1^4,2^4,3^4,4)$ with $10\times\frac12$ points and such that \[T=W\cap H_1\cap H_2\cap H_3\cap H_4,\] where the $H_i$ are hyperplanes of the projective space $\PP(1^4,2^4,3^4,4)$.\end{mainthm}
The Main Theorem is proved in Section \ref{sec!hellextend}, using a projection--unprojection construction for $T$ and for $W^6$. See \cite{PR} and \cite{Ki} for details on projection--unprojection methods.

We call the Fano $6$-fold $W$ a \emph{key variety}, because it contains several interesting varieties as appropriate weighted complete intersections. We already know how to recover the curve $D$ and the K3 surface $T$ from $W^6$, and we can also construct hyperelliptic surfaces of general type.
\begin{cor}\label{cor!covering} There is a $15$-parameter family of hyperelliptic surfaces $Y$ of general type with $p_g=1$, $q=0$, $K^2=2$ and no torsion, each of which is a complete intersection of type $(1,1,1,2)$ in a Fano $6$-fold $W\subset\PP(1^4,2^4,3^4,4)$ with $10\times\frac12$ points.
\end{cor}
Each surface $Y$ contains the genus $3$ hyperelliptic curve $D\in|K_Y|$ as the unique hyperplane section of weight $1$. These surfaces were constructed by Catanese and Debarre in \cite{CDE}, but our method has the advantage of being more widely applicable to other situations.

In Section \ref{sec!generaltype} we use the key variety method to construct a new family of Godeaux surfaces with torsion $\ZZ/2$.
\begin{thm} There is an $8$-parameter family of Godeaux surfaces $X$ with torsion $\ZZ/2$, where each $X$ is obtained as a $\ZZ/2$-quotient of some hyper\-elliptic surface $Y$ constructed in Corollary \ref{cor!covering}.
\end{thm}
The Godeaux surfaces are surfaces of general type with $p_g=0$, $K^2=1$, and their torsion group is cyclic of order $\le5$. The components of the moduli space with torsion $\ZZ/5$, $\ZZ/4$, $\ZZ/3$ were constructed in \cite{R1}, and in each case the moduli space is irreducible, unirational, and $8$-dimensional, which is the expected dimension. The first simply connected example appeared in \cite{B}, and recently another simply connected Godeaux surface was constructed in \cite{LP} using $\QQ$-Gorenstein smoothing theory. It is expected that the key variety method will give an irreducible $8$-dimensional component of the moduli space of Godeaux surfaces with algebraic fundamental group $\ZZ/2$, although we have not proved that here.

\subsection*{Acknowledgements} I would like to thank Miles Reid for introducing me to this problem, which forms part of my University of Warwick PhD thesis \cite{mythesis}. This research was partially supported by the World Class University program through the National Research Foundation of Korea funded by the Ministry of Education, Science and Technology (R33-2008-000-10101-0).

\section{Graded rings over hyperelliptic curves}\label{sec!hellcurve}
In this section we review hyperelliptic curves and their graded rings, using the double covering of $\PP^1$. This is well known material covered in section $4$ of \cite{Aq}. We include this section for completeness since we will generalise to hyperelliptic K3 surfaces in Section~\ref{sec!K3surface}.

We consider the case when $D$ is a hyperelliptic curve of genus $3$. Then the canonical linear system $|K_D|$ defines a double covering of $\PP^1$ embedded as a plane conic, branched in $8$ points $Q_1,\dots,Q_8$. The corresponding ramification (or Weierstrass) points on $D$ are labelled $P_1,\dots,P_8$. The double covering $\pi\colon D\to\PP^1$ determines and is determined by the $g^1_2$: a free linear system of dimension $1$ and degree $2$. Moreover, we have $2P_i\sim g^1_2$ and $K_D\sim2g^1_2$. There is a natural hyperelliptic involution $h$ on $D$ which swaps the two sheets of the double covering, and $\pi$ is the quotient map of this involution.

Choose generators $s_1$, $s_2$ of $H^0(D,g^1_2)$. These are coordinates on $\PP^1$, and there is a polynomial $F_8(s_1,s_2)$ whose vanishing determines the branch locus $Q_1+\dots+Q_8$ of $\pi$. The double covering is $D_8\subset\PP(1,1,4)$, defined by the equation $w^2=F_8(s_1,s_2)$. By considering rational functions on $D$, we have \[4g^1_2\sim P_1+\dots+P_8,\] or more generally, \[P_1+\dots+P_a + (8-a)g^1_2\sim P_{a+1}+\dots+P_8 + 4g^1_2.\]

We write $B_1=P_1+\dots+P_a$, $B_2=P_{a+1}+\dots+P_8$, then since the $B_i$ are effective Cartier divisors on $D$ we can choose constant sections \[u\colon\Oh_D\to\Oh_D(B_1),\quad  v\colon\Oh_D\to\Oh_D(B_2).\] Now $u^2, uv, v^2$ are sections of $ag^1_2$, $4g^1_2$, $(8-a)g^1_2$ respectively, so we have two relations $u^2=f(s_1,s_2)$, $v^2=g(s_1,s_2)$ and the identity $w=uv$. Here $f(s_1,s_2)$ is a homogeneous function of degree $a$ on $\PP^1$ with zeros at $Q_1,\dots,Q_a$, similarly $g(s_1,s_2)$, so that $F=fg$.

Clearly every $h$-invariant divisor class can be written in the form \[A\sim P_1+\dots+P_a+bg^1_2\sim P_{a+1}+\dots+P_8+(a+b-4)g^1_2.\] For such a divisor $A$, the graded ring \[R(D,A)=\bigoplus_{n\ge0} H^0(D,\Oh_D(nA))\] can be studied relative to the base $\PP^1$ via the double covering $\pi$. We quote the following proposition from \cite{Aq} for $D$ of genus $3$, although the proposition and subsequent graded ring calculations work for any genus with only minor alterations.
\begin{prop}\label{prop!hellcurve} Let $D$ be a hyperelliptic curve of genus $3$ with Weierstrass points $P_1,\dots,P_8$, and write $\pi\colon D\to\PP^1$ for the natural quotient by the hyperelliptic involution $h$. Then
\begin{enumerate}
\item[(1)] $\pi_*\Oh_D=\Oh_{\PP^1}\oplus\Oh_{\PP^1}(-4);$
\item[(2)] $\pi_*\Oh_D(g^1_2)=\Oh_{\PP^1}(1)\oplus\Oh_{\PP^1}(-3);$
\item[(3)] $\pi_*\Oh_D(P_1+\dots+P_a)=\Oh_{\PP^1}u\oplus\Oh_{\PP^1}(a-4)v;$
\end{enumerate}
where in each case the first summand is invariant under $h$ and the second is anti-invariant.\end{prop}
\begin{rmk}\rm Note that in case (1) the direct image sheaf is a sheaf of $\Oh_{\PP^1}$-algebras, where the multiplication \[\Oh_{\PP^1}(-4)\otimes\Oh_{\PP^1}(-4)\to\Oh_{\PP^1}\] is defined via $w^2=F(s_1,s_2)$.
\end{rmk}
\subsection{The ineffective theta characteristic}\label{sec!curveconstruction}
Consider the ineffective theta characteristic \[A\sim P_1+\dots+P_4-g^1_2\sim P_5+\dots+P_8-g^1_2\] on $D$. An ineffective theta characteristic is a divisor class on $D$ such that $h^0(D,A)=0$ and $2A\sim K_D$. Using Proposition \ref{prop!hellcurve}, we see that $R(D,A)$ is generated by monomials in $s_1, s_2, u, v$:
\[\begin{array}{c|c|c}
n & H^0(D,\Oh_D(nA)) & H^0(\PP^1,\pi_*\Oh_D(nA))\\
\hline
0 & 1 & 1 \\
1 & \phi & \phi \\
2 & y_1,\ y_2,\ y_3 & s_1^2,\ s_1s_2,\ s_2^2\\
3 & z_1,\ z_2,\ z_3,\ z_4 & s_1u,\ s_2u,\ s_1v,\ s_2v\\
4 & t & uv \\
\end{array}\]
The relations between these generators are either of the trivial monomial kind, or derived from \[u^2=f_4(s_1,s_2),\ v^2=g_4(s_1,s_2).\] For example, it is clear that $z_1^2=s_1^2u^2=y_1f(y_1,y_2,y_3)$, where $f(y_1,y_2,y_3)$ is a \emph{rendering} of $f_4(s_1,s_2)$ in the quadratic monomials $s_1^2$, $s_1s_2$, $s_2^2$. In fact, we can present all the equations as
\[\rank\left(\renewcommand{\arraystretch}{1.2}\begin{array}{cc|cc}
y_1&y_2&z_1&z_3\\
y_2&y_3&z_2&z_4\\
\hline
z_1&z_2&f_2&t\\
z_3&z_4&t&g_2
\end{array}\right)\le1,\] where $f_2$ and $g_2$ are quadrics in $y_1$, $y_2$, $y_3$. Taking $\Proj R(D,A)$ gives \[D\subset\PP(2^3,3^4,4),\] and the double covering of $\PP^1$ is the conic defined by the first $2\times2$ minor of the matrix.

\section{Graded rings over hyperelliptic K3 surfaces}\label{sec!K3surface}
In this section we generalise the methods of Section \ref{sec!hellcurve} to hyperelliptic K3 surfaces, and in \ref{sec!K3construction} we construct a hyperelliptic K3 surface $T$ that extends the hyperelliptic curve $D$ of \ref{sec!curveconstruction}. Section \ref{sec!Cstar} gives alternative descriptions of $D$ and $T$, while \ref{sec!K3proj} describes a projection construction for $T$ which will be used in the proof of the Main Theorem \ref{thm!hell}.

A hyperelliptic K3 surface $T$ is a K3 surface together with a complete linear system $L$ such that $|L|$ contains an irreducible hyperelliptic curve $D$ of arithmetic genus $g=h^0(T,\Oh_T(L))-1$. Then $L$ determines a $2$-to-$1$ map $\pi\colon T\to F$ where $F$ is a surface of degree $g-1$ in $\PP^g$. The branch locus of the double covering is some divisor in $|{-}2K_F|$. For further details on the hyperelliptic dichotomy for K3 surfaces see \cite{SD}. Del Pezzo classified the possibilities for $F$ as rational scrolls or the Veronese surface. Since both have very simple explicit descriptions, we can analyse graded rings over any hyperelliptic K3 surface by calculating relative to the base $F$. For brevity we treat only the case $g=3$, but more general examples are contained in \cite{mythesis}.

We assume that $F=Q_2\subset\PP^3$ is a quadric of rank $4$ and the double cover $\pi\colon T\to F$ is branched in a curve $C$ of bidegree $(4,4)$, which splits into two components $C_1+C_2$ of bidegree $(3,1)$ and $(1,3)$ respectively. The components of the branch curve intersect one another transversally in $10$ points which are nodes of $T$. This is the hyperelliptic degeneration of the symmetric determinantal quartic K3 surface of \cite{part1}. As usual there is a hyperelliptic involution $h\colon T\to T$ exchanging the two sheets of the double cover, and $\pi$ is the quotient map of $h$. Let $H_1$, $H_2$ be the generators of $\Pic Q$, then we omit $\pi^*$ to write $\pi^*H_i=H_i$ on $T$, and $\pi^*C_i=2D_i$.

Let $s_1, s_2$ be generators of $H^0(T,H_1)$, similarly $t_1, t_2$ for $H^0(T,H_2)$. Then there is an equation $F_{4,4}(s_1,s_2,t_1,t_2)$ defining the branch curve $C$ on $Q$. This equation factors as $F=f_{3,1}(s_i,t_i)g_{1,3}(s_i,t_i)$, which determines the splitting $C=C_1+C_2$. The double cover $T$ is given by $w^2=F$, and we have $2D_1\sim3H_1+H_2$ and $2D_2\sim H_1+3H_2$ on $T$. Considering the rational function $w/(t_1^2s_1^2)$ on $T$, we find \[2(H_1+H_2)\sim D_1+D_2.\]

By analogy with the hyperelliptic curves of section \ref{sec!hellcurve} we write down graded rings \[R(T,A)=\bigoplus_{n\ge0}H^0(T,\Oh_T(nA))\] where $A$ is a divisor class which is invariant under $h$. Any such $A$ can be written in the form \[A\sim D_1+n_1H_1+n_2H_2\sim D_2+(n_1+1)H_1+(n_2-1)H_2.\]The following proposition is a natural extension of Proposition \ref{prop!hellcurve}, which allows us to describe $R(T,A)$ relative to $R(Q,\pi_*A)$.

\begin{prop}\label{prop!hellK3} Let $T$ be a hyperelliptic K3 surface double covering of the rank $4$ quadric $Q\subset\PP^3$, with ramification properties as described above. Choose constant sections $u\colon\Oh_T\to\Oh_T(D_1)$ and $v\colon\Oh_T\to\Oh_T(D_2)$ for the components $D_i$ of the ramification curve. Clearly we have $u^2=f_{3,1}(s_i,t_i)$, $uv=w$ and $v^2=g_{1,3}(s_i,t_i)$, where $F=fg$. Moreover,
\begin{itemize}
\item[(1)] $\pi_*\Oh_T=\Oh_Q\oplus\Oh_Q(-2,-2)$;
\item[(2)] $\pi_*\Oh_T(H_1)=\Oh_Q(1,0)\oplus\Oh_Q(-1,-2)$;
\item[(3)] $\pi_*\Oh_T(D_1)=\Oh_Qu\oplus\Oh_Q(1,-1)v$;
\end{itemize}
with similar results for $H_2$, $D_2$ respectively.
\end{prop}
\begin{rmk}\rm Once again we note the $\Oh_Q$-algebra structure on $\pi_*\Oh_T$. The multiplication map \[\Oh_Q(-2,-2)\otimes\Oh_Q(-2,-2)\to\Oh_Q\] is defined via the equation $w^2=F_{4,4}(s_1,s_2,t_1,t_2)$.
\end{rmk}

\subsection{Construction of the K3 surface $T$}\label{sec!K3construction}
Write $A\sim D_1-H_1\sim D_2-H_2$, which is an $h$-invariant divisor class on $T$, satisfying $H^0(T,\Oh_T(A))=0$ and $\Oh_T(A)^{[2]}=\pi^*\Oh_Q(1)$.  Note that $A$ is the analogue of the ineffective theta characteristic in Section \ref{sec!curveconstruction}. We can describe the ring $R(T,A)$ using Proposition \ref{prop!hellK3}. The generators for $R(T,A)$ are:
\[\renewcommand{\arraystretch}{1.2}\begin{array}{c|c|c}
n & H^0(T,\Oh_T(nA) & H^0(Q,\pi_*\Oh_T(nA))\\
\hline
0 & 1 & 1 \\
1 & 0 & 0 \\
2 & y_1,\ y_2,\ y_3,\ y_4 & s_1t_1,\ s_2t_1,\ s_1t_2,\ s_2t_2\\
3 & z_1,\ z_2,\ z_3,\ z_4 & t_1u,\ t_2u,\ s_1v,\ s_2v\\
4 & t & uv=w \\
\end{array}\]
The relations are again mostly trivial monomial relations, together with those derived from $u^2=f_{3,1}$ and $v^2=g_{1,3}$. Some are slightly more difficult to write down than others, for example,
\[z_1t=t_1u^2v=t_1vf_{3,1}=s_1vq_{2,2}+s_2vq^\prime_{2,2}=z_3q_1(y_i)+z_4q_1^\prime(y_i),\]
where $q_1$ and $q_1^\prime$ are suitable quadrics rendered in $y_1$,\dots,$y_4$.
The trick here is to make $f_{3,1}$ bihomogeneous by incorporating the factor $t_1$ into $f$ and simultaneously taking out the excess in $s_1$, $s_2$. Clearly we can not expect $f$ to be divisible by $s_1$ or by $s_2$, and so we have a choice of ways to break up $f$ into quadrics. Fortunately, this choice is arbitrary, as any discrepancy is accounted for by the rank condition (\ref{eq!wtedsegre}) below. We present all the relations of $R(T,A)$ as follows
\begin{equation}\label{eq!wtedsegre}\rank\left(\renewcommand{\arraystretch}{1.2}\begin{array}{cc|c}
y_1&y_2&z_1\\
y_3&y_4&z_2\\
\hline
z_3&z_4&t
\end{array}\right)\le1,\end{equation}
\begin{align*}
  z_1^2&=t_1^2f_{3,1}&  z_3^2&=s_1^2g_{1,3}\\
z_1z_2&=t_1t_2f_{3,1}&  z_3z_4&=s_1s_2g_{1,3}\\
  z_2^2&=t_2^2f_{3,1}&  z_4^2&=s_2^2g_{1,3}\\
  z_1t&=q_1z_3+q_1^\prime z_4 &  z_3t&=q_3z_1+q_3^\prime z_2 \\
  z_2t&=q_2z_3+q_2^\prime z_4 &  z_4t&=q_4z_1+q_4^\prime z_2
\end{align*}
\[t^2=F(y_i),\]
where for example $z_1^2=t_1^2f_{3,1}$ means we render the bihomogeneous expression $t_1^2f_{3,1}$ in the variables $y_1,\dots,y_4$. Then $\Proj R(T,A)$ gives us the K3 surface \[T\subset\PP(2^4,3^4,4),\] which has $10\times\frac12(1,1)$ points. Note that the curve $D$ of Section \ref{sec!curveconstruction} is obtained by taking a hyperplane section of weight $2$ in $T$, avoiding the $\frac12$ points.

\subsection{Alternative descriptions of hyperelliptic varieties}\label{sec!Cstar}

We can consider the curve $D$ as a codimension $2$ complete intersection inside a weighted homogeneous variety as follows: let $X$ be the second Veronese embedding of $\PP^3$ with coordinates $s_1, s_2, u, v,$ and take the affine cone $\Ce X\subset\AA^{10}$ over $X$. Aside from the obvious $\CC^\times$-action on $\Ce X$ there are many other possibilities, and we choose a weighted $\CC^\times$-action with weights $(1,1,2,2)$. Then the quotient $Y=\Ce X\git_{\!1}\CC^\times$ of $\Ce X$ is contained in $\PP(2^3,3^4,4^3)$, and is defined by the equations 
\[\rank\left(\renewcommand{\arraystretch}{1.2}\begin{array}{cc|cc}
y_1&y_2&z_1&z_3\\
y_2&y_3&z_2&z_4\\
\hline
z_1&z_2&x_1&t\\
z_3&z_4&t&x_2
\end{array}\right)\le1.\]
The hyperelliptic curve $D$ is simply the codimension $2$ complete intersection $x_1=f_2$, $x_2=g_2$ inside $Y$.

Similarly, $T$ is a codimension $2$ complete intersection in the weighted homogeneous space we now describe. Consider the following $(\CC^\times)^2$-action on $\CC^2\times\CC^2\times\CC^2$:
\begin{align*}\la\colon (s_1,s_2,t_1,t_2,u,v)&\mapsto(\la^2s_1,\la^2s_2,t_1,t_2,\la^3u,\la v)\\
\mu\colon (s_1,s_2,t_1,t_2,u,v)&\mapsto(s_1,s_2,\mu^2t_1,\mu^2t_2,\mu u,\mu^3v).\end{align*}
The $4$-dimensional quotient $Z=(\CC^2\times\CC^2\times\CC^2)\git_{\!(1,1)}(\CC^\times)^2$ is embedded in $\PP(2^4,3^4,4)$ by the determinantal equations (\ref{eq!wtedsegre}). The surface $T$ is the complete intersection $u^2=f_{3,1}\in(6,2)$, $v^2=g_{1,3}\in(2,6)$ in $Z$.

\subsection{Projection of $T$ to a complete intersection}\label{sec!K3proj}

Consider the following projection map on del Pezzo surfaces: let $Q\subset\PP^3$ be a quadric of rank $4$ and blow up a point $P$ on $Q$ to obtain the del Pezzo surface $B$. Then contract the two $(-1)$-curves on $B$ arising from the rulings of $Q$ to get $\PP^2$. Now, suppose we have a curve $C$ on $Q$ of type $(4,4)$ which splits as $C=C_1+C_2$ where $C_1\in(3,1)$, $C_2\in(1,3)$ so that $C$ has $10$ nodes. If the centre of projection $P$ is chosen to be one of these nodes then the two components $C_1$, $C_2$ are projected to nodal plane cubics, and the image of $P$ is the line $L$ through these two nodes.

Now suppose we have a hyperelliptic $K3$ surface $T$ which is a double cover of $Q$ branched in $C$. The classical projection $Q\dasharrow\PP^2$ lifts to the double cover as illustrated by the diagram below:
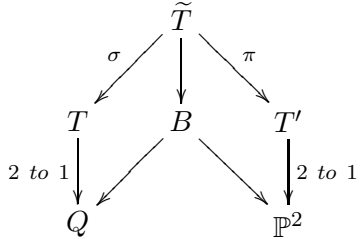
\[\xymatrix{&\widetilde{T}\ar[dl]_\sigma\ar[d]\ar[dr]^\pi&\\
            T\ar[d]_{2\ to\ 1}&B\ar[dl]\ar[dr]&T^\prime\ar[d]^{2\ to\ 1}\\
            Q&&\PP^2}\]
where $\sigma\colon\widetilde{T}\to T$ is the blowup of $P$ in $T$ and we write $E\cong\PP^1$ for the exceptional divisor. The image $T^\prime$ of the projection is a double cover of $\PP^2$ branched over the two nodal cubics. The centre of projection $P$ in $T$ is projected to a rational curve of arithmetic genus $2$ double covering $L$ away from the two nodes, and branched over the residual intersection with $C$.

Now this diagram can also be recast as a projection--unprojection operation in the sense of \cite{Ki}, \cite{PR}. Start from $T\subset\PP(2^4,3^4,4)$ with $10\times\frac12$ points and polarising divisor $A$, as described in section \ref{sec!K3construction}. Choose a $\frac12$ point $P$ in $T$, and write $\sigma\colon\widetilde{T}\to T$ for the $(1,1)$-weighted blowup of $P$, whose exceptional curve is $E$. Then the projected surface $T^\prime_{6,6}\subset\PP(2,2,2,3,3)$ is calculated as 
\[T^\prime=\Proj R(\widetilde T,\sigma^*A-\textstyle{\frac12}E),\]
where certain functions on $\widetilde T$ are eliminated by the projection since they do not vanish appropriately along $E$. Geometrically, the surface $T_{6,6}^\prime$ is a double covering of $\PP^2$ branched in the two nodal cubics defined by the equations of degree $6$.

This projection to $T^\prime$ can be expressed explicitly as an operation in commutative algebra. Assume the centre of projection is a $\frac12$ point at the coordinate point $P_{y_4}$, with local coordinates $z_3,z_4$. Then adjusting the notation of section \ref{sec!K3construction} slightly, write down the matrix relations
\begin{equation}\label{eq!projhellK3matrix}
\rank\renewcommand{\arraystretch}{1.2}\begin{pmatrix}
y_2&f&z_1\\
g&y_4&z_3\\
z_2&z_4&t\end{pmatrix}\le1,\end{equation}
where we reserve the right to choose $f$, $g$ later. These equations are a subset of those for $T\subset\PP(2^4,3^4,4)$ after a trivial change of coordinates. The remaining equations for $T$ are completely determined by 
\begin{align*}z_1^2&=L_1y_2^2+L_2y_2f+L_3f^2\\z_2^2&=M_1y_2^2+M_2y_2g+M_3g^2\end{align*}
where a priori $L_i$, $M_i$ are linear in $y_1,\dots,y_4$. Indeed, the equations we have written down so far are sufficient to determine the two components of the branch curve, and their defining equations $f_{3,1}$ and $g_{1,3}$. We can fill in the remaining equations of $T$ using the procedure outlined in Section \ref{sec!K3construction}.

Since we fixed a $\frac12$ point at $P_{y_4}$, the last equation for $T$ can be written as \[t^2=a_2(y_1,y_3)y_4^2+b_3(y_1,y_2,y_3)y_4+c_4(y_1,y_2,y_3).\] Now the tangent cone to $P$ must factorise because the branch curve $C$ splits into two components, so we can choose coordinates \[f=y_1+\al y_3,\quad g=\be y_1+y_3\] so that $a=y_1y_3$. This in turn forces $L_3=y_1$, $M_3=y_3$ so that modulo the minors of matrix (\ref{eq!projhellK3matrix}), the equations involving $z_1^2$ and $z_2^2$ take the form
\begin{equation}\begin{split}\label{eqn!hellprojected}
z_1^2&=L_1(y_1,y_2,y_3)y_2^2+l_4y_2fg+y_1f^2\\
z_2^2&=M_1(y_1,y_2,y_3)y_2^2+m_4y_2fg+y_3g^2,\end{split}\end{equation}
where $L_1$, $M_1$ do not involve $y_4$ and $l_4$, $m_4$ are scalars. The image of the exceptional curve $E$ is defined by $y_2=0$.

We are finally in a position to describe the projection centred at $P_{y_4}$ in terms of explicit equations. The local coordinates near $P$ are $z_3$, $z_4$ so we expect the projection to eliminate these variables along with $y_4$ (see \cite{Ki}, example 9.13). In fact the projection also eliminates $t$, and we are left with equations (\ref{eqn!hellprojected}) defining a complete intersection \[T^\prime_{6,6}\subset\PP(2,2,2,3,3).\] This is the hyperelliptic degeneration of the \emph{totally tangent conic} configuration of \cite{part1}.

\section{Extending hyperelliptic graded rings}\label{sec!hellextend}
In this section we consider extensions of the hyperelliptic K3 surface $T$ constructed in section \ref{sec!K3construction}, and prove the Main Theorem \ref{thm!hell}. As with the symmetric determinantal extensions of \cite{part1}, the most convenient way to extend the K3 surface $T$ is by using the projection construction of Section \ref{sec!K3proj}. We start from
\[\PP^1\xrightarrow{\,\fie\,}T^\prime_{6,6}\subset\PP(2,2,2,3,3),\] where $T^\prime$ is a double covering of $\PP(2,2,2)$ branched in two nodal cubics. The image of $\fie$ is a curve of arithmetic genus $2$, which is a double cover of the line joining the two nodes. Constructing $\fie$ and $T^\prime_{6,6}$ is equivalent to constructing $T$ itself, so we prove the theorem by extending $\fie$ and $T^\prime$.

We assume that $\fie$ is a double cover of the line $(y_2=0)\subset\PP(2,2,2)$ branched over the points $\fie(1,0)$ and $\fie(0,1)$. Then for general $T^\prime$ the map $\fie$ is
\[\fie\colon\PP^1\to\PP(2,2,2,3,3)\]
\begin{equation}\label{eqn!hellfiemap}(u,v)\mapsto\left(u^2,0,v^2,u(u^2+\al v^2),v(\be u^2+v^2)\right).\end{equation}
Rendering $\fie^*(z_i^2)$ in terms of $y_1$, $y_3$ we see that the image of $\fie$ is defined by the equations
\begin{align}
C_1\colon z_1^2&=y_1(y_1+\al y_3)^2\label{eqn!hell1}\\
C_2\colon z_2^2&=y_3(\be y_1+y_3)^2\label{eqn!hell2}\\
y_2&=0.\label{eqn!hell3}
\end{align}
To define $T^\prime\subset\PP(2,2,2,3,3)$ we must choose two appropriate combinations of weight $6$ in equations (\ref{eqn!hell1}--\ref{eqn!hell3}). Note that if we want the branch curves to be nondegenerate then we should ensure that both equations for $T^\prime$ involve $y_2$ nontrivially. Moreover, after incorporating $y_2$ into the equations we should check that there are still two bona fide nodes on the branch locus at $(-\al,0,1)$ and $(1,0,-\be)$. So, calculating the tangent cone to each curve at these points forces the equations of $T^\prime$ to take the form
\begin{equation}\label{eqn!hellTprime}\begin{split}C_1&+l_1Q_1+l_2Q_2+l_3Q_3+l_4Q_4\\
C_2&+m_1Q_1+m_2Q_2+m_3Q_3+m_4Q_4,
\end{split}\end{equation}
where $\al$, $\be$, $l_i$, $m_i$ are scalar parameters and
\[Q_1=(y_1+\al y_3)y_2^2,\quad Q_2=y_2^3,\quad Q_3=(\be y_1+y_3)y_2^2,\]
\[Q_4=(y_1+\al y_3)(\be y_1+y_3)y_2.\]

\begin{pfmainthm} The proof follows a similar approach to the main result of \cite{part1} and it is informative to compare the two at each stage. We explicitly extend the projected image $T^\prime_{6,6}\subset\PP(2,2,2,3,3)$ to a Fano $6$-fold $W_{6,6}^\prime\subset\PP(1^4,2^3,3^2)$ containing the image of $\PP^5$ under some map $\Fie$. Define $\fie\colon\PP^1\to\PP(2,2,2,3,3)$ as in (\ref{eqn!hellfiemap}) and write $\fie_0\colon\PP^1\to\PP(2,2,2)$ for the map \[\fie_0^*(y_1)=u^2,\quad\fie_0^*(y_2)=0,\quad\fie_0^*(y_3)=v^2.\] Then writing $u$, $v$, $a$, $b$, $c$, $d$ for the coordinates on $\PP^5$, up to automorphisms of $\PP^5$ and $\PP(1^4,2^3)$ the general extension of $\fie_0$ to $\Fie_0\colon\PP^5\to\PP(1^4,2^3)$ is
\[\Fie_0^*(a)=a,\quad\Fie_0^*(b)=b,\quad\Fie_0^*(c)=c,\quad\Fie_0^*(d)=d,\]
\begin{equation}\label{def!hellFie0}\setlength\arraycolsep{2pt}\renewcommand{\arraystretch}{1.3}\begin{array}{rlrr}\Fie_0^*(y_1)&=u^2&&+2av,\\\Fie_0^*(y_2)&=0\phantom{^2}+&bu&+\phantom{2}cv,\\\Fie_0^*(y_3)&=v^2+&2du&\end{array}\end{equation}
We prove that there is a unique map $\Fie\colon\PP^5\to\PP(1^4,2^4,3^2)$ which is a lift of $\Fie_0$ and which extends $T^\prime_{6,6}$ to $W^\prime_{6,6}$.

Write $M$, $R$, $S$ for the coordinate rings of $\PP^5$, $\PP(1^4,2^3)$ and $\PP(1^4,2^3,3^2)$ respectively. By equation (\ref{def!hellFie0}), the map $\Fie_0^*$ induces a graded $R$-module structure on $M$ with generators $1$, $u$, $v$ and $uv$. Similarly $\Fie^*$ makes $M$ into a graded $S$-module with the same generators. The presentation of $M$ as a module over $R$ is
\[0\leftarrow M\xleftarrow{(1,u,v,uv)}R\oplus2R(-1)\oplus R(-2)\xleftarrow{A}R(-2)\oplus2R(-3)\oplus R(-4)\]
where $A$ is the matrix
\begin{equation}\label{matrix!hellA}\renewcommand{\arraystretch}{1.2}\begin{pmatrix}
-y_2&by_1&cy_3&-2cdy_1+4ady_2-2aby_3\\
b&-y_2&-2cd&cy_3\\
c&-2ab&-y_2&by_1\\
0&c&b&-y_2\end{pmatrix}.\end{equation}
Since $\Fie$ is a lift of $\fie$ we assume that the general forms of $\Fie^*(z_i)$ are
\begin{align*}
\Fie^*(z_1)&=u^3+\al uv^2+s_1u^2+s_2uv+s_3v^2+s_4u+s_5v\\
\Fie^*(z_2)&=\be u^2v+v^3+t_1u^2+t_2uv+t_3v^2+t_4u+t_5v
\end{align*}
where the $s_i(a,b,c,d)$, $t_i(a,b,c,d)$ are homogeneous polynomials of degree $1$ or $2$ as appropriate. Then using the $R$-module structure of $M$ we can write
\begin{equation}\begin{split}\label{def!hellFie}
\Fie^*(z_1)&=u(f+s_4)+s_2uv+s_5v\\
\Fie^*(z_2)&=v(g+t_5)+t_2uv+t_4u
\end{split}\end{equation}
where \[f=y_1+\al y_3,\quad g=\be y_1+y_3.\]
We have used coordinate changes $z_1\mapsto z_1+s_1y_1$ and similar to absorb the values of $s_1$, $s_3$, $t_1$, $t_3$ into $z_1$, $z_2$. The following theorem shows that there are unique values of $s_i$, $t_i$ for $i=2,4,5$ for which there are equations extending (\ref{eqn!hell1}), (\ref{eqn!hell2}). As a corollary, we prove that for these unique values of $s_i$, $t_i$, there are extensions of equations $Q_1,\dots,Q_4$.
\begin{thm}\label{thm!hellextend}\begin{enumerate}\item[(I)]The kernel of $\Fie^*\colon S\to M$ contains equations extending (\ref{eqn!hell1}), (\ref{eqn!hell2}) of the form \begin{align*}z_1^2-y_1f^2&\in R+Rz_1+Rz_2,\\z_2^2-y_3g^2&\in R+Rz_1+Rz_2\end{align*} if and only if 
\begin{align*}s_2=(1-\al\be)a,\quad s_4=\be a^2+\al^2d^2,\quad s_5=\al(\al\be-1)ad,\\
t_2=(1-\al\be)d,\quad t_4=\be(\al\be-1)ad,\quad t_5=\be^2a^2+\al d^2.\end{align*}
\item[(II)] Given part $(I)$, the equations are
\begin{align}\begin{split}\label{eq!hellkerFie1}z_1^2-y_1(f+s_4)^2&=-4(f+s_4)s_2ay_3 - 4s_2s_5dy_1 + s_2^2y_1y_3 + s_5^2y_3\\&\quad+ 2(1-\al\be)a^2(3dz_1-az_2) - 2\al a(f+s_4)z_2\end{split}\\
\begin{split}\label{eq!hellkerFie2}z_2^2-y_3(g+t_5)^2&=-4(g+t_5)t_2dy_1 - 4t_2t_4ay_3 + t_2^2y_1y_3 + t_4^2y_1\\
&\quad+ 2(1-\al\be)d^2(3az_2-dz_1) - 2\be d(g+t_5)z_1\end{split}
\end{align}
\end{enumerate}\end{thm}
\begin{cor}\label{cor!hellextend}The kernel of $\Fie^*$ also contains (nontrivial) equations extending $Q_i$ for $i=1,\dots,4$ of the form
\[fy_2^2,\ y_2^3,\ gy_2^2,\ fgy_2\in R+Rz_1+Rz_2\] respectively.
\end{cor}
\begin{pf} The ``if'' part of the theorem is proved by evaluating equations (\ref{eq!hellkerFie1}), (\ref{eq!hellkerFie2}) under $\Fie^*$ with $s_i$, $t_i$ taking the values stated in the theorem. The remainder of the proof is for the ``only if'' part.

Using the graded module structure of $k[u,v]$ over $k[y_1,y_2,y_3]$ via $\fie_0^*$ we write
\begin{align*}\fie^*(z_1)&=(y_1+\al y_3)u\\\fie^*(z_2)&=(\be y_1+y_3)v.\end{align*} Then squaring either of these expressions and rendering $u^2$, $v^2$ as $y_1$, $y_3$ gives equations (\ref{eqn!hell1}), (\ref{eqn!hell2}) immediately.
We attempt to do the same rendering calculation for the extended map $\Fie^*$, using
\begin{align*}u^2&=\Fie^*(y_1)-2av\\v^2&=\Fie^*(y_3)-2du.\end{align*} We can eliminate all terms involving $u^2$ or $v^2$ from $\Fie^*(z_i^2)$ to obtain
\begin{align*}\Fie^*\left(z_1^2-y_1(f+s_4)^2+4(f+s_4)s_2ay_3 + 4s_2s_5dy_1 - s_2^2y_1y_3 - s_5^2y_3\right)&\equiv0\\\Fie^*\left(z_2^2-y_3(g+t_5)^2+4(g+t_5)t_2dy_1 + 4t_2t_4ay_3 - t_2^2y_1y_3 - t_4^2y_1\right)&\equiv0\end{align*} modulo $(a,b,c,d)M$. The residual parts to these congruences are
\begin{align*}K&=K_uu+K_vv+K_{uv}uv,\\L&=L_uu+L_vv+L_{uv}uv\end{align*}
respectively, where
\begin{equation}\begin{split}\label{eq!hellresK}K_u&=8(f+s_4)s_2ad-2s_5^2d-2s_2^2dy_1+2s_2s_5y_3\\K_v&=-2(f+s_4)^2a+8s_2s_5ad+2(f+s_4)s_2y_1-2s_2^2ay_3\\K_{uv}&=2(f+s_4)s_5+4s_2^2ad\end{split}\end{equation}
and
\begin{equation}\begin{split}\label{eq!hellresL}L_u&=-2(g+t_5)^2d + 8t_2t_4ad + 2(g+t_5)t_2y_3 - 2t_2^2dy_1\\L_v&=8(g+t_5)t_2ad - 2t_4^2a - 2t_2^2ay_3 + 2t_2t_4y_1\\L_{uv}&=2(g+t_5)t_4+4t_2^2ad.\end{split}\end{equation}

Now $K$, $L$ are homogeneous expressions of degree $6$ in $(a,b,c,d)M$, and we prove that if they are to be contained in the submodule $R+Rz_1+Rz_2\subset M$ then $s_i$, $t_i$ must take the values stated in the theorem. From the definition of $\Fie^*(z_i)$ in (\ref{def!hellFie}), the submodule $R+Rz_1+Rz_2$ is the image of the composite
\[M\xleftarrow{(1,u,v,uv)}R\oplus2R(-1)\oplus R(-2)\xleftarrow{B}R\oplus2R(-3)\oplus R(-2)\oplus2R(-3)\oplus R(-4)\]
where $B$ is the matrix
\[\renewcommand{\arraystretch}{1.2}\left(\begin{array}{ccc|cccc}1&0&0&-y_2&by_1&cy_3&-2cdy_1+4ady_2-2aby_3\\
0&f+s_4&t_4&b&-y_2&-2cd&cy_3\\
0&s_5&g+t_5&c&-2ab&-y_2&by_1\\
0&s_2&t_2&0&c&b&-y_2\end{array}\right).\]
The first $3$ columns of $B$ are the generators $1$, $z_1$, $z_2$ and the last $4$ columns are the matrix $A$ from (\ref{matrix!hellA}), which is mapped to $0$ under the composite.

We seek vectors $\xi,\eta\in R\oplus2R(-3)\oplus R(-2)\oplus2R(-3)\oplus R(-4)$ such that
\begin{equation}\begin{split}\label{eq!findxieta}K=\begin{pmatrix}1,&u,&v,&uv\end{pmatrix}B\xi,\\
L=\begin{pmatrix}1,&u,&v,&uv\end{pmatrix}B\eta.\end{split}\end{equation} In order to solve for $\xi$, $\eta$ and consequently fix the values of $s_i$, $t_i$ we stratify $K,L$ according to degree in $y_1,y_2,y_3$. In other words, write
\begin{align*}K&=K^{(0)}+K^{(1)}+K^{(2)}\\L&=L^{(0)}+L^{(1)}+L^{(2)}\end{align*} where $K^{(i)}$, $L^{(i)}$ have degree $i$ in $y_1,y_2,y_3$ and similarly we write \begin{align*}\xi&=\xi^{(0)}+\xi^{(1)}\\\eta&=\eta^{(0)}+\eta^{(1)}.\end{align*}

We begin with $K^{(2)}$, which is calculated from (\ref{eq!hellresK}) as \[K^{(2)}=2f(y_1s_2-fa)v.\] We must find $\xi^{(1)}$ such that \begin{equation}\label{eq!calcK2}K^{(2)}=\begin{pmatrix}1,&u,&v,&uv\end{pmatrix}B\xi^{(1)} + \text{ lower order terms.}\end{equation} Comparing coefficients of $y_1^2$ and $y_3^2$, the only solution is \[\xi_3^{(1)}=\frac2\be(s_2-a)y_1-2\al^2ay_3,\]with the other $\xi_i^{(1)}=0$. Then the coefficient of $y_1y_3$ in (\ref{eq!calcK2}) dictates that \[s_2=(1-\al\be)a\]and therefore $\xi_3^{(1)}=-2\al af$. An exactly similar calculation with $L^{(2)}$ and $\eta_2^{(1)}$ yields \[t_2=(1-\al\be)d\]and $\eta_2^{(1)}=-2\be dg$.

Proceeding to the calculation for $K^{(1)}$, we must solve
\begin{equation}\label{eq!calcK1}K^{(1)}-\xi_3^{(1)}(t_4u+t_5v+t_2uv)=\begin{pmatrix}1,&u,&v,&uv\end{pmatrix}B\xi^{(0)} + \text{ lower order terms}\end{equation} where the term involving $\xi_3^{(1)}$ is necessary to account for the lower order terms from equation (\ref{eq!calcK2}).
Now examining the coefficient of $uv$ in (\ref{eq!calcK1}), we obtain
\[2f(s_5+\al at_2)=s_2\xi^{(0)}_2+t_2\xi^{(0)}_3.\] However, $\xi^{(0)}$ has degree $0$ in $y_i$ by construction, so the left hand side must be identically $0$. Hence \[s_5=-\al at_2\] and by considering the coefficient of $uv$ in $L^{(1)}$ we find \[t_4=-\be ds_2.\]

Comparing coefficients of $u$ and $v$ in equation (\ref{eq!calcK1}) we obtain
\begin{align*}6(1-\al\be)a^2df&=(f+s_4)\xi_2^{(0)}+t_4\xi_3^{(0)}+\text{ lower order terms}\\
2a(-s_4(f+\al g)+\al ft_5-s_2^2y_3)&=s_5\xi_2^{(0)}+(g+t_5)\xi_3^{(0)}+\text{ lower order terms}.
\end{align*}Since $\xi^{(0)}$ has degree $0$ in $y_i$ we must have $\xi_2^{(0)}=6(1-\al\be)a^2d$. Moreover the coefficient of $v$ must be divisible by $g$, which is equivalent to \begin{equation}\label{eq!simult1}\al t_5-s_4=-\be(1-\al\be)a^2.\end{equation} By considering the coefficients of $u$, $v$ in $L^{(1)}$ in the same way we get $\eta_3^{(0)}=6(1-\al\be)ad^2$ and a further restriction on $s_4$, $t_5$:
\begin{equation}\label{eq!simult2}t_5-\be s_4=\al(1-\al\be)d^2.\end{equation}

Solving equations (\ref{eq!simult1}), (\ref{eq!simult2}) simultaneously forces
\begin{align*}s_4&=\be a^2+\al^2d^2\\t_5&=\be^2a^2+\al d^2,\end{align*}
which in turn means that
\begin{align*}\xi_3^{(0)}&=-2(1-\al\be)a^3-2\al as_4\\
\eta_2^{(0)}&=-2(1-\al\be)d^3-2\be dt_5.
\end{align*}
We can finally write out $\xi$ and $\eta$ in full
\begin{align*}\xi_2&=6(1-\al\be)a^2d&\eta_2&=-2\be d(g+t_5)-2(1-\al\be)d^3\\
\xi_3&=-2\al a(f+s_4)-2(1-\al\be)a^3&\eta_3&=6(1-\al\be)ad^2,
\end{align*}where the other $\xi_i=\eta_i=0$. It is necessary to check that $\xi$ and $\eta$ actually solve equations (\ref{eq!findxieta}) when all the lower order terms are replaced, which can be verified directly.

The extended equations (\ref{eq!hellkerFie1}), (\ref{eq!hellkerFie2}) are obtained by writing out the vectors $\xi$, $\eta$ in terms of the generators of $R+Rz_1+Rz_2$
\begin{align*}\begin{split}z_1^2-y_1(f+s_4)^2&=-4(f+s_4)s_2ay_3 - 4s_2s_5dy_1 + s_2^2y_1y_3 + s_5^2y_3\\&\quad+\xi_2z_1+\xi_3z_2\\
z_2^2-y_3(g+t_5)^2&=-4(g+t_5)t_2dy_1 - 4t_2t_4ay_3 + t_2^2y_1y_3 + t_4^2y_1\\
&\quad+\eta_2z_1+\eta_3z_2.\end{split}\end{align*}
This concludes the proof of theorem (\ref{thm!hellextend}).
\end{pf}
\begin{pfcor} First observe that the fourth column of $B$ is equivalent to $y_2=bu+cv$. Thus the extension of $Q_1$ is calculated by expressing $fy_2(bu+cv)$ in terms of the other columns of $B$. We have to find $\nu$ such that
\[(f+s_4)y_2(bu+cv)=\begin{pmatrix}1,&u,&v,&uv\end{pmatrix}B\nu.\] The solution to this linear algebra problem is\begin{align*}
\nu_2&=2by_2+2(\be ab-cd)c&\nu_3&=2(\al b^2+c^2)a\\
\nu_4&=-\be as_2y_2-2ac(g+t_5)+2(cd-\be ab)s_5&\nu_5&=b(f+s_4)-\be abs_2\\
\nu_6&=-c(f+s_4)-\be acs_2+2bs_5&\nu_7&=2bs_2,\end{align*}
where $\nu_1=y_2\nu_4-by_1\nu_5-cy_3\nu_6-(-2cdy_1+4ady_2-2aby_3)\nu_7$ uses the first column of $B$ to remove any excess terms. Thus the equation extending $Q_1$ is
\[\widetilde Q_1\colon(f+s_4)y_2^2=\nu_1+\nu_2z_1+\nu_3z_2.\]

Similar calculations give the equations extending $Q_2$, $Q_3$, $Q_4$ for which we list the corresponding vectors below. The equation extending $Q_2$ is
\[\widetilde Q_2\colon y_2^3=\nu_1+\nu_2z_1+\nu_3z_2,\] where
\begin{align*}
\nu_1&=y_2\nu_4-by_1\nu_5-cy_3\nu_6-(-2cdy_1+4ady_2-2aby_3)\nu_7\\
\nu_2&=\frac 2{\al\be - 1}(b^2 + \be c^2)b\\
\nu_3&=\frac 2{\al\be - 1}(\al b^2 + c^2)c\\
\nu_4&=\frac{2}{1 - \al\be}\left(b^2(f + s_4)+c^2(g + t_5)\right) + (\be ac + \al bd)y_2 + 2(2 - \al\be)abcd\\
\nu_5&=-by_2 + 2c^2d + (\be ac + \al bd)b\\
\nu_6&=-cy_2 + 2ab^2 + (\be ac + \al bd)c\\
\nu_7&=-2bc.
\end{align*}
The equation extending $Q_3$ is
\[\widetilde Q_3\colon (g+t_5)y_2^2=\nu_1+\nu_2z_1+\nu_3z_2,\] where
\[\nu_1=y_2\nu_4-by_1\nu_5-cy_3\nu_6-(-2cdy_1+4ady_2-2aby_3)\nu_7\]
\begin{align*}
\nu_2&=2(b^2 + \be c^2)d&\nu_3&=2cy_2 - 2(ab - \al cd)b\\
\nu_4&=-\al dt_2y_2 - 2bd(f + s_4) + 2(ab - \al cd)t_4&
\nu_5&=-b(g + t_5) - \al bdt_2 + 2ct_4\\
\nu_6&=c(g + t_5) - \al cdt_2&\nu_7&=2ct_2.
\end{align*}
Finally, equation $Q_4$ is extended by
\[\widetilde Q_4\colon (f+s_4)(g+t_5)y_2=\nu_1+\nu_2z_1+\nu_3z_2\]
where
\[\nu_1=y_2\nu_4-by_1\nu_5-cy_3\nu_6-(-2cdy_1+4ady_2-2aby_3)\nu_7\]
\begin{align*}
\nu_2&=b(g + t_5) + ct_4 - t_2y_2&\nu_3&=c(f + s_4) + bs_5 - s_2y_2\\
\nu_4&=-s_5t_4&\nu_5&=-t_2(f + s_4) - s_2t_4\\
\nu_6&=-s_2(g + t_5) - s_5t_2&\nu_7&=-2s_2t_2.
\end{align*}
This completes the proof of the corollary.
\end{pfcor}
Given Theorem \ref{thm!hellextend} and its Corollary, we can prove that there is a unique hyperelliptic Fano $6$-fold $W^\prime_{6,6}\subset\PP(1^4,2^3,3^2)$ extending any given projected hyperelliptic K3 surface $T^\prime_{6,6}$. Simply take the combination of equations (\ref{eq!hellkerFie1}), (\ref{eq!hellkerFie2}) and $\widetilde Q_i$ corresponding to the choice (\ref{eqn!hellTprime}) made in the definition of $T^\prime_{6,6}$. This proves the Main Theorem \ref{thm!hell}.
\end{pfmainthm}

\section{Godeaux surfaces with torsion $\ZZ/2$}\label{sec!generaltype}
In this section we describe a family of hyperelliptic surfaces $Y$ of general type with $p_g=1$, $q=0$, and $K^2=2$. Each surface $Y$ is a complete intersection inside a key variety $W$ constructed in Main Theorem \ref{thm!hell}. We then find an appropriate subfamily of surfaces that are Galois \'etale $\ZZ/2$-coverings of Godeaux surfaces. In particular, we give an explicit description of the fixed point free $\ZZ/2$-action on $Y$, which is the restriction of an appropriate $\ZZ/2$-action on the key variety $W$.

\subsection{Coverings of Godeaux surfaces}
Let $X$ be the canonical model of a surface of general type with $p_g=0$, $K^2=1$. We call $X$ a Godeaux surface, and we assume that the torsion subgroup $\Tors X\subset\Pic X$ has order $2$. Write $\sigma$ for the generator of $\Tors X$, and consider the Galois \'etale double covering $f\colon Y\to X$ induced by $\sigma$. The covering surface is constructed by taking $Y=\Proj R(Y,K_X,\sigma)$, or written out in full
\[Y=\Proj\bigoplus_{n\ge0}\left(H^0(X,nK_X)\oplus H^0(X,nK_X+\sigma)\right).\]
The surface $Y$ is the canonical model of a surface of general type with $p_g=1$, $q=0$, $K^2=2$, and the extra $\ZZ/2$-grading on the ring $R(Y,K_X,\sigma)$ determines a fixed point free $\ZZ/2$ group action on $Y$, where the first summand is invariant and the second is anti-invariant. The quotient by this group action is the map $f$.

Moreover, an analysis of $R(Y,K_X,\sigma)$ reveals that the canonical curve of $Y$ must be hyperelliptic.
\begin{lemma}\label{lem!hyperelliptic} If $Y$ is the unramified double covering of a Godeaux surface with torsion $\ZZ/2$, then the canonical curve section $D$ in $|K_Y|$ is hyperelliptic.
\end{lemma}
This lemma was also proved in \cite{CDE}, using a monomial counting proof. We use a Hilbert series approach which has some advantages, and yields slightly more information about the group action for use later.
\begin{pf} Define the bigraded Hilbert series of the ring $R(Y,K_X,\sigma)$ by
\[P_Y(t,e)=\sum_{n\ge0}\left(h^0(X,nK_X)t^n+h^0(X,nK_X+\sigma)t^ne\right),\]
where $t$ keeps track of the degree, and $e$ keeps track of the eigenspace, so that $e^2=1$. Then using the Riemann--Roch theorem,
\[P_Y(t,e)=1+et+2t^2+2t^2e+4t^3+4t^3e+\dots\]
which can be written as the rational function
\[P_Y(t,e)=\frac{1+(e-1)t^4+(-2e-2)t^5+(-4e-6)t^6+(7e+8)t^8+\dots}{(1-et)(1-t^2)(1-et^2)^2(1-t^3)^2(1-et^3)^2}.\]
Using well-known Hilbert series properties would normally indicate that $Y$ is a subvariety of $\PP(1,2^3,3^4)$. However, the first nontrivial coefficient in the numerator is not negative, due to the bigrading. Thus we must introduce an extra generator of degree $4$ in the negative eigenspace, dividing $P_Y(t,e)$ by $(1-et^4)$ so that the numerator becomes \[1-t^4+(-2e-2)t^5+(-4e-6)t^6+\dots.\] The extra $-t^4$ term in the numerator suggests that it is necessary to introduce a relation in degree $4$, which does not eliminate the new generator of degree $4$. Thus the canonical curve section of $Y$ must be hyperelliptic, and is given by $D\subset\PP(2^4,3^4,4)$.
\end{pf}
Now, we can construct hyperelliptic surfaces $Y$ of general type using the key variety of Main Theorem \ref{thm!hell}.
\begin{thm}\label{thm!coverings}There is a $15$-parameter family of hyperelliptic surfaces $Y$ of general type with $p_g=1$, $q=0$, $K^2=2$ and no torsion, each of which is a complete intersection of type $(1,1,1,2)$ in a Fano $6$-fold $W\subset\PP(1^4,2^4,3^4,4)$ with $10\times\frac12$ points.
\end{thm}
The proof is identical to that of \cite{part1}, Theorem 4.1, and we obtain the canonical model of $Y$ using this construction. In order to find Godeaux surfaces with $\ZZ/2$-torsion, we must determine which hyperelliptic surfaces $Y$ have an appropriate $\ZZ/2$-action. To do this we study the $\ZZ/2$-action on the hyperelliptic curve $D$, and extend it to the key variety $W$.

\subsection{The canonical curve}

Now let $f\colon Y\to X$ be the \'etale double cover of a Godeaux surface $X$ and suppose $D$ is a nonsingular curve in $|K_Y|$, similarly $C$ in $|K_X+\sigma|$. By Lemma \ref{lem!hyperelliptic}, $D$ is a hyperelliptic curve and $C$ has genus $2$ so is automatically hyperelliptic too. Let $\pi_D\colon D\to Q\cong\PP^1$ denote the quotient map of the hyperelliptic involution on $D$, similarly $\pi_C\colon C\to\PP^1$. Since $D$ is an unramified double cover of $C$ via $f|_D$, this induces a double covering of $\Image\pi_C=\PP^1$ by $Q$.

We get the following picture:
\[\xymatrix{
E\ar[r]&D\ar[r]^{f|_D}\ar[d]^{\pi_D}&C\ar[d]^{\pi_C}\\
&Q\ar[r]&\PP^1}\]
There is a fixed point free involution on the curve $D$ induced by the unramified double cover $f|_D$, which we call the Godeaux involution. We use the same notation for the Godeaux involution and the torsion element $\sigma\in\Pic X$. The Weierstrass points of $D$ must be invariant under $\sigma$, so there is a natural division of these eight points into two sets $\{P_1,\dots,P_4\}$ and $\{P_5,\dots,P_8\}$, which are interchanged by $\sigma$.

Now consider the ineffective theta characteristic $A_D=K_Y|_D$ on $D$ which is determined by the surface $Y$. The divisor class $A_D$ is invariant under both $\sigma$ and the hyperelliptic involution, so a priori the only possibilities are 
\begin{equation}\label{eqn!AD}\begin{split}A_D&\sim P_1+P_2+P_3+P_4-g^1_2\sim P_5+P_6+P_7+P_8-g^1_2,\\
A_D&\sim P_1+P_3+P_5+P_7-g^1_2\sim P_2+P_4+P_6+P_8-g^1_2.\end{split}\end{equation}
The difference between these is that the former is only $\sigma$-invariant as a divisor class, whereas the latter is an $\sigma$-invariant divisor.

Now, we have already constructed the graded ring $R(D,A_D)$ in section \ref{sec!curveconstruction}. Furthermore by the adjunction formula, it is clear that $2g^1_2\sim 2A_D\sim K_D$. However, these two divisor classes $A_D$ and $g^1_2$ are distinct, because the $g^1_2$ is effective whereas $A_D$ is ineffective. Thus we have a $2$-torsion class \[\tau=A_D-g^1_2\] on $D$, which corresponds to a genus $5$ unramified double cover $E$ of $D$, where\[E=\Proj R(D,A_D,\tau)=\Proj\bigoplus_{n\ge0}\left(H^0(D,nA_D)\oplus H^0(D,nA_D+\tau)\right).\]
We outline the procedure to construct the bigraded ring $R(D,A_D,\tau)$. Using the notation of section \ref{sec!hellcurve}, write $s_1,s_2$ for the sections of the $g^1_2$ and \[u\colon\Oh_D\to\Oh_D(P_1+\dots+P_4),\quad v\colon\Oh_D\to\Oh_D(P_5+\dots+P_8).\] We can very quickly write down generators and relations for $R(D,A_D,\tau)$:
\[\renewcommand{\arraystretch}{1.1}\begin{array}{c|c|c}
n&H^0(D,nA_D)&H^0(D,nA_D+\tau)\\
\hline
0&k&\phi\\
1&\phi&s_1,s_2\\
2&s_1^2,s_1s_2,s_2^2&u,v\\
3&\dots&\dots
\end{array}\] Thus $E$ is a complete intersection \[E_{4,4}\subset\PP(1,1,2,2),\] defined by equations $u^2=f_4(s_1,s_2)$ and $v^2=g_4(s_1,s_2)$. The polynomials $f$ and $g$ are functions on $\PP^1$ whose vanishing determines the splitting of the Weierstrass points of $D$ into two sets of four.

The curve $E$ comes bundled at no extra cost with the fixed point free involution $\tau\colon E\to E$ associated to the torsion $\tau$ of $D$. We recover the restricted algebra $R(D,K_Y|_D)$ of section \ref{sec!curveconstruction} by taking the $\tau$-invariant subring of $R(D,A_D,\tau)$:\[R(D,A_D)=R(D,A_D,\tau)^{\left<\tau\right>}.\] For future reference, we write out the action of $\tau$ on $E$ using the eigenspace table above \[s_1\mapsto -s_1,\quad s_2\mapsto -s_2,\quad u\mapsto -u,\quad v\mapsto -v.\]

Now, we claim that the covering curve $E$ completely determines the Godeaux involution $\sigma$ on $D$. First observe that $D$ is a quotient of $E$, and that this covering curve only exists because $D$ is the curve section of $|K_Y|$. Thus $\sigma$ lifts to the curve $E$ and should be compatible with the involution $\tau$ on $E$, so that $\sigma^2=1$ or $\tau$ on $E$.

\begin{prop}\label{prop!sigmacurve} The action of $\sigma$ on $E$ is given by
\[s_1\mapsto is_1,\quad s_2\mapsto -is_2,\quad u\mapsto iv,\quad v\mapsto iu,\] so that $\sigma^2=\tau$ and the group $\left<\sigma,\tau\right>$ is isomorphic to $\ZZ/4$. Moreover, the polarising divisor of $D$ is 
\[A_D\sim P_1+P_2+P_3+P_4-g^1_2\sim P_5+P_6+P_7+P_8-g^1_2.\]
\end{prop}
\begin{pf}
The Hilbert series of Lemma \ref{lem!hyperelliptic} gives the eigenspace decomposition of $\sigma$ on $D$, which we must abide by. In particular, $R(D,A_D)$ should have only one invariant generator in degree $2$, and the generator in degree $4$ should be anti-invariant. This forces $\sigma^2=\tau$, so that the group $\left<\sigma,\tau\right>$ acting on $E$ is $\ZZ/4$ rather than $\ZZ/2\oplus\ZZ/2$.

Now, there are two possibilities for $\sigma$ depending on the representation of $A_D$ chosen from equation (\ref{eqn!AD}). The correct choice is  \[s_1\mapsto is_1,\quad s_2\mapsto -is_2,\quad u\mapsto iv,\quad v\mapsto iu,\] which corresponds to the first possibility in (\ref{eqn!AD}).
Indeed, the alternative \[s_1\mapsto is_1,\quad s_2\mapsto -is_2,\quad u\mapsto iu,\quad v\mapsto iv,\] is obliged to have two fixed points on $D$ at the coordinate points $P_{y_1}$ and $P_{y_3}$, and so can not possibly be the Godeaux involution.

Hence we have \[A_D\sim P_1+P_2+P_3+P_4-g^1_2\sim P_5+P_6+P_7+P_8-g^1_2,\] with corresponding action on $R(D,A_D)$ given by
\[\renewcommand{\arraystretch}{1.2}\begin{array}{ccc}
\rank\left(\begin{array}{cc|cc}
y_1&y_2&z_1&z_3\\
y_2&y_3&z_2&z_4\\
\hline
z_1&z_2&f_2&t\\
z_3&z_4&t&g_2
\end{array}\right)\le1&
\mapsto&
\rank\left(\begin{array}{cc|cc}
-y_1&y_2&-z_3&-z_1\\
y_2&-y_3&z_4&z_2\\
\hline
-z_3&z_4&-g_2&-t\\
-z_1&z_2&-t&-f_2
\end{array}\right)\le1,\end{array}\]
where
\begin{align*}
f_2&=\al_1y_1^2+\al_2y_1y_2+\al_3y_1y_3+\al_4y_2^2+\al_5y_2y_3+\al_6y_3^2\\
g_2&=-\al_1y_1^2+\al_2y_1y_2-\al_3y_1y_3-\al_4y_2^2+\al_5y_2y_3-\al_6y_3^2\end{align*}
and the involution has no fixed points as long as $\al_1$ and $\al_6$ are not zero. This proves the proposition.
\end{pf}
\subsection{Involution on the K3 surface}\label{sec!invK3}
Moving one step up the tower, we lift the involution $\sigma$ on the canonical curve section $D$ to the hyperelliptic K3 surface $T\subset\PP(2^4,3^4,4)$, which contains $D$ as a quadric section.

The whole argument becomes quite transparent when viewed in terms of commutative algebra. The graded ring $R(T,A_T)$ is described explicitly in section \ref{sec!K3construction}, and after eliminating one of the generators in degree $2$, we obtain $R(D,A_D)$. Now if $D$ is the unramified double covering of a Godeaux curve $C$ with its involution $\sigma\colon D\to D$ from Proposition \ref{prop!sigmacurve}, then:
\begin{prop}\label{prop!sigmaK3} There is at least one K3 surface $T$ containing the curve $D$ such that the involution $\sigma$ on $D$ has a unique lift to $T$. Moreover, such a lift $\sigma\colon T\to T$ has four fixed points which are $\frac12$ points of $T$. We call $\sigma$ the Godeaux involution on $T$.
\end{prop}
\begin{rmk}\rm This is surprising because we are looking for a fixed point free involution on the covering surface $Y$, so it would be reasonable to expect that the involution on the K3 surface is free.
\end{rmk}
\begin{pf}\emph{Step (1) Determining the character of $\sigma$.} Temporarily choose coordinates on $T$ so that $D=T\cap(y_4=0)$, where $y_4$ must be semi-invariant under any putative involution. Then the determinantal equations (\ref{eq!wtedsegre}), which partially define $T$, take the general form
\[\renewcommand{\arraystretch}{1.2}
\rank\left(\begin{array}{ccc}
y_1+\al y_4&y_2+\be y_4&z_1\\
y_2+\gam y_4&y_3+\delta y_4&z_2\\
z_3&z_4&t
\end{array}\right)\le1,\] where $\al,\be,\gam,\delta$ are scalars. Now if $\sigma$ lifts to $T$, then our choice of coordinates means that the action of $\sigma$ on $T$ is predetermined by $\sigma|_D$ from Proposition \ref{prop!sigmacurve}, excepting the new variable $y_4$. Since $T$ is a double covering of a quadric $Q\subset\PP^3$ of rank $4$, the determinantal equations force $\al=\delta$ and $\be=-\gam$. Thus $y_4$ is anti-invariant, and the signature of $\sigma$ on $Q$ is $(1,3)$. We can recalibrate the coordinate system so that the determinantal equations and involution on $T$ are
\[\renewcommand{\arraystretch}{1.2}\begin{array}{ccc}
\rank\left(\begin{array}{cc|c}
y_1&y_2&z_1\\
y_3&y_4&z_2\\
\hline
z_3&z_4&t
\end{array}\right)\le1
&\mapsto&
\rank\left(\begin{array}{cc|c}
-y_1&y_3&-z_3\\
y_2&-y_4&z_4\\
\hline
-z_1&z_2&-t
\end{array}\right)\le1,\end{array}\]
where the original curve $D$ is obtained from $T$ by taking the anti-invariant quadric section $y_2=y_3$.\vspace{0.3cm}

\noindent\emph{Step (2) Fixed points of $\sigma$.} First observe that the involution on $T$ swaps the two branch curves, and also swaps the sheets of the hyperelliptic double covering $\pi\colon T\to Q$. Thus any fixed points of $\sigma$ lie on both components of the branch curve, and so must be $\frac12$ points of $T$.

For a $\frac12$ point of $T$ to be fixed under $\sigma$, one of two things must happen:
\[y_1=y_4=y_2-y_3=0\text{, or }y_2+y_3=0.\]
The only case we need to worry about is when $y_2+y_3=0$ since the other case reduces to the curve $D$, on which $\sigma$ is fixed point free by hypothesis. To ensure $C_1$ and $C_2$ are interchanged under $\sigma$, their equations are of the form
\begin{equation}
\begin{split}\label{eqn!branch}
f_{3,1}&=\al_1s_1^3t_1+\al_2s_1^2s_2t_1+\al_3s_1s_2^2t_1+\al_4s_2^3t_1\\
&\qquad+\be_1s_1^3t_2+\be_2s_1^2s_2t_2+\be_3s_1s_2^2t_2+\be_4s_2^3t_2,\\
g_{1,3}&={-}\al_1s_1t_1^3+\al_2s_1t_1^2t_2-\al_3s_1t_1t_2^2+\al_4s_1t_2^3\\
&\qquad+\be_1s_2t_1^3-\be_2s_2t_1^2t_2+\be_3s_2t_1t_2^2-\be_4s_2t_2^3,
\end{split}
\end{equation}
on $Q=\PP^1\times\PP^1$. Note that there is more than one possible choice of $f_{3,1}$, $g_{1,3}$ for which $T\cap(y_2=y_3)=D$, so we can not claim that $T$ is unique in the statement of the proposition. We can assume $y_2=1$, $y_3={-}1$, then for a general choice of branch curve there are four fixed points on $T$. These are $(\la,1,-1,-1/\la,0,0,0,0,0)$, where $\la$ is a root of the quartic equation derived from evaluating equations (\ref{eqn!branch})
\[\al_1\la^4+(\al_2-\be_1)\la^3+(\al_3-\be_2)\la^2+(\al_4-\be_3)\la-\be_4,\]
which proves the proposition.
\end{pf}
\subsection{Involution on the Fano 6-fold}\label{sec!invFano}
We extend the involution on the K3 surface $T$ to the Fano $6$-fold $W$ constructed in Main Theorem \ref{thm!hell}. To do this we use a $\ZZ/2$-equivariant form of the projection--unprojection construction described in Section \ref{sec!K3proj}.

We begin with a $\ZZ/2$-equivariant unprojection construction for the K3 surface $T$ with a Godeaux involution. Recall from Proposition \ref{prop!sigmaK3} that if $T$ has a Godeaux involution $\sigma$, then $\sigma$ has four fixed points, each of which is a $\frac12$ point. Let $P$ be one of these fixed points, and project from $P$, to obtain \[\PP^1\xrightarrow{\fie}T_{6,6}^\prime\subset\PP(2,2,2,3,3).\] Here $T^\prime$ is a double cover of $\PP^2$ branched in two nodal cubics, and $P$ is projected to the image of $\fie$, a double covering of the line through the two nodes. There is an induced involution on $T^\prime$, which swaps the two branch cubics and leaves the image of $\fie$ invariant. We call this a $\ZZ/2$\emph{-equivariant projection}.

Examining the equations of $T^\prime$ as described in Section \ref{sec!K3proj}, the determinantal equations for $T$ are
\[\rank\renewcommand{\arraystretch}{1.2}\begin{pmatrix}
y_2&f&z_1\\
g&y_4&z_3\\
z_2&z_4&t\end{pmatrix}\le1,\]
where $f=y_1+\al y_3$, $g=\al y_1+y_3$ because the two branch curves are interchanged by $\sigma$. Proposition \ref{prop!sigmaK3} fixes the involution on $T$ as
\begin{equation}\begin{split}\label{def!sigmaTprime}f\mapsto g&,\quad y_2\mapsto -y_2,\quad g\mapsto f,\quad z_1\mapsto-z_2,\quad z_2\mapsto-z_1,\\
&y_4\mapsto-y_4,\quad z_3\mapsto z_4,\quad z_4\mapsto z_3,\quad t\mapsto-t\end{split}\end{equation}
and we note that this implies $\sigma(y_1)=y_3$, $\sigma(y_3)=y_1$. Hence referring to Section \ref{sec!K3proj} the equations of $T^\prime_{6,6}\subset\PP(2,2,2,3,3)$ must be of the form
\begin{equation}\begin{split}\label{eq!sigmaTprime}
z_1^2&=y_1f^2+y_2^2(l_1f+l_2y_2+l_3g)+l_4y_2fg\\
z_2^2&=y_3g^2+y_2^2(l_3f-l_2y_2+l_1g)-l_4y_2fg,
\end{split}\end{equation} where $l_i$ are scalars. The remaining equations of $T$ can be calculated from those of $T^\prime$ using the unprojection procedure outlined in Section \ref{sec!K3proj}.

There are three isolated fixed points on $T^\prime$ when $z_1=z_2=y_1+y_3=0$, which correspond to three of the $9\times\frac12$ points as expected. Further, $T^\prime$ has two fixed points on the unprojection divisor which arise from the fact that the centre of projection $P$ was itself a fixed point. Indeed, suppose we have a local orbifold chart for a neighbourhood of the $\frac12$ point $P$ in $T$. This is the quotient of $\CC^2$ by $\ZZ/2$ acting by ${-}1$ on both coordinates. Then writing $u,v$ for the coordinates on $\CC^2$, $\sigma$ lifts to the chart as \[u\mapsto-iv,\quad v\mapsto-iu\]by equation (\ref{def!sigmaTprime}). The Kawamata $(1,1)$ weighted blowup at $P$ introduces the ratio $(u:v)$ as the exceptional $\PP^1$, which is then embedded in $T^\prime$ by $\fie$. Thus the induced action of $\sigma$ on the image of $\fie$ inside $T^\prime$ has two fixed points at $\fie(1,1)$ and $\fie(-1,1)$.

It is important to note that the $\ZZ/2$-equivariant projection--unprojection construction for $T$ relies on the choice of $\frac12$ point $P$. As such we can no longer assume there is a canonical choice of curve $D\subset T$ defined by setting $f=g$, as we did in the proof of Proposition \ref{prop!sigmaK3}. The choice of covering curve $D$ is made by taking any anti-invariant quadric section of $T$ which avoids the $10\times\frac12$ points. Indeed, the quadric $f=g$ contains the point $P$ and so is no longer a valid choice.

Now, we claim that the involution on $T$ can be extended to the Fano $6$-fold $W$ at the top of the tower.
\begin{prop}\label{prop!sigmaFano} Suppose $T\subset\PP(2^4,3^4,4)$ is a K3 surface with $10\times\frac12$ points and $\sigma\colon T\to T$ is a Godeaux involution lifted from some quadric section $D\subset T$. Then there is a lift of $\sigma$ to the unique Fano $6$-fold $W\subset\PP(1^4,2^4,3^4,4)$ extending $T$ which was constructed in Main Theorem \ref{thm!hell}. Moreover the fixed locus of the involution $\sigma\colon W\to W$ consists of four isolated $\frac12$ points.
\end{prop}
\begin{pf} Project from one of the fixed $\frac12$ points on $T$ to get \[\fie\colon\PP^1\to T^\prime_{6,6}\subset\PP(2^3,3^2).\] Following the extension procedure outlined in the proof of Main Theorem \ref{thm!hell}, the extended map \[\Fie\colon\PP^5\to W^\prime_{6,6}\subset\PP(1^4,2^3,3^2)\] must be \[\Fie\colon(a,b,c,d,u,v)\mapsto(a,b,c,d,u^2+2av,bu+cv,v^2+2du,f_1,f_2),\]
where
\begin{align*}f_1 &= u\left(f + \al(a^2 + \al d^2)\right) + (1-\al^2)auv + \al(\al^2-1)adv,\\
      f_2 &= v\left(g + \al(\al a^2+d^2)\right) + (1 - \al^2)duv + \al(\al^2-1)adu.
\end{align*}
To make $\Fie$ compatible with the lift of $\sigma\colon T\to T$ defined by equation (\ref{def!sigmaTprime}), the action on $\PP^5$ must be
\[u\mapsto -v,\quad v\mapsto-u,\quad a\mapsto -d,\quad b\mapsto c,\quad c\mapsto b,\quad d\mapsto -a.\]
Thus $\Fie$ is $\sigma$-equivariant, so the equations defining the image of $\Fie$ are invariant and consequently $W^\prime\subset\PP(1^4,2^3,3^2)$ can be chosen to be invariant. Alter\-natively, a direct calculation following the proof of Main Theorem \ref{thm!hell} demonstrates explicitly that the equations of the image of $\Fie$ are invariant. Hence by $\ZZ/2$-equivariant unprojection, the involution lifts to the $6$-fold $W$.

Now outside the image of $\Fie$, there are just three isolated points on $W^\prime$ that are fixed under $\sigma$. These are the same $\frac12$ points that were fixed under $\sigma|_{T^\prime}$. On the image of $\Fie$ itself there are two copies of $\PP^2\subset\PP^5$ whose image under $\Fie$ are fixed by $\sigma$. These are defined by \[\PP^5\cap(u=v,a=d,b=-c),\quad\PP^5\cap(u=-v,a=-d,b=c),\] and they are the analogue of the two fixed points on $\fie(\PP^1)\subset T^\prime$. Fortunately these nonisolated fixed loci are contracted to the centre of projection $P$ on $W$, so that $\sigma$ fixes just four isolated $\frac12$ points there. This proves the proposition.
\end{pf}
\section{Godeaux surfaces with torsion $\ZZ/2$}

Given a hyperelliptic tower $D\subset T\subset W$ where $W$ is the unique Fano $6$-fold extending the K3 surface $T$, suppose the curve $D$ is a double cover of a Godeaux curve $C$. Now, suppose further that the tower is constructed so that the Godeaux involution $\sigma$ on $D$ lifts to $T$ and subsequently $W$ as described in Propositions \ref{prop!sigmaK3} and \ref{prop!sigmaFano}. Write $A$ for the hyperplane class on $W$ so that $\Oh_W(A)=\Oh_W(1)$, and $-K_W=4A$. Then $\sigma$ induces a $\ZZ\oplus\ZZ/2$-bigrading on the ring $R(W,A)$ according to eigenspace:
\[\begin{array}{c|c|c}
n&H^0(W,nA)^+&H^0(W,nA)^-\\
\hline
1&a-d, b+c & a+d, b-c\\
2&y_1+y_3&y_1-y_3,y_2,y_4\\
3&z_1-z_2,z_3+z_4&z_1+z_2,z_3-z_4\\
4&&t
\end{array}\]
Now by Theorem \ref{thm!coverings}, we can construct a surface $Y$ of general type with $p_g=1$, $K^2=2$ as a complete intersection inside $W$ as long as $Y$ avoids the $\frac12$ points of $W$, which is an open condition. Referring to the above table and the eigenspace decomposition on $Y$ given by lemma \ref{lem!hyperelliptic}, if we $Y$ to be a complete intersection of type $(1^+,1^+,1^-,2^-)$ inside $W$ then $\sigma|_Y$ will be the fixed point free Godeaux involution. Hence we have:
\begin{thm}There is an $8$ parameter family of Godeaux surfaces with $\ZZ/2$-torsion.
\end{thm}
The parameter count is a matter of calculating the moduli of $W$ using Main Theorem \ref{thm!hell}, Section \ref{sec!invFano} and then counting the number of free parameters involved in choosing the complete intersection $(1^+,1^+,1^-,2^-)$.

\bigskip\noindent
Stephen Coughlan, \\
Math. Dept., Sogang University, \\
Sinsu-dong, Mapo-gu, \\
Seoul, Korea \\
e-mail: stephencoughlan21@gmail.com \\

\end{document}